\documentclass[12pt,a4paper]{amsart}
\usepackage[utf8]{inputenc}
\usepackage[english]{babel}
\usepackage[foot]{amsaddr}
\usepackage{amsmath,amsfonts,amssymb,amsthm}
\usepackage{mathtools} 
\usepackage{dsfont} 
\usepackage[percent]{overpic}
\usepackage[cal=euler]{mathalfa}
\usepackage[dvipsnames]{xcolor}
\usepackage{hyperref}
\hypersetup{
    colorlinks = true,
    linkcolor = {RedOrange},
    citecolor = {ForestGreen},
}

\newtheorem{theorem}{Theorem}
\newtheorem{proposition}[theorem]{Proposition}
\newtheorem{corollary}[theorem]{Corollary}
\newtheorem{lemma}[theorem]{Lemma}

\theoremstyle{definition}

\newtheorem{example}[theorem]{Example}

\theoremstyle{remark}
\newtheorem{remark}[theorem]{Remark}

\DeclareMathOperator{\Mod}{Mod}
\DeclareMathOperator{\Stab}{Stab}

\DeclareMathOperator{\tw}{tw}
\DeclareMathOperator{\lo}{o}
\DeclareMathOperator{\bO}{O}
\DeclareMathOperator{\Homeo}{Homeo}
\DeclareMathOperator{\GL}{GL}
\DeclareMathOperator{\msp}{m}

\newcommand{\Th}{\mathrm{Th}}
\newcommand{\WP}{\mathrm{WP}}

\newcommand{\BB}{\mathbb{B}}

\newcommand{\ZZ}{\mathbb{Z}}
\newcommand{\QQ}{\mathbb{Q}}
\newcommand{\RR}{\mathbb{R}}
\newcommand{\bS}{\mathbb{S}}

\newcommand{\PP}{\mathbb{P}}

\newcommand{\cB}{\mathcal{B}}

\newcommand{\cM}{\mathcal{M}}
\newcommand{\cML}{\mathcal{ML}}

\newcommand{\cP}{\mathcal{P}}

\newcommand{\cS}{\mathcal{S}}
\newcommand{\cT}{\mathcal{T}}

\begin{document}
\title[Length statistics of random multicurves on closed hyperbolic surfaces]{Length statistics of random multicurves \\ on closed hyperbolic surfaces}
\author{Mingkun Liu}
\address{Université de Paris and Sorbonne Université, CNRS, IMJ-PRG, F-75006 Paris, France}
\email{mingkun.liu@imj-prg.fr}

\begin{abstract}
    A multicurve can be decomposed into components.
    The ratios of the length of each component to the total length give a hint of the shape of the multicurve.
    In this paper, we determine the distribution of these ratios of a random multicurve with a given topological type on a closed hyperbolic surface, using the methods of Margulis' thesis and Mirzakhani's equidistribution theorem for horospheres.
    This distribution admits a polynomial density whose coefficients can be expressed explicitly in terms of intersection numbers of psi-classes on the Deligne--Mumford compactification of the moduli space of complex curves, and in particular it does not depend on the hyperbolic metric.
    This result generalizes prior work of M.~Mirzakhani in the case of random pants decompositions. Results very close to ours were obtained independently and simultaneously by F.~Arana-Herrera.
\end{abstract}
\maketitle

\section{Introduction}
Let $X$ be a connected closed oriented complete hyperbolic surface of genus $g \geq 2$.
An \emph{ordered multi-geodesic} on $X$ is a finite ordered list $\gamma = (m_1\gamma_1,\dots,m_k\gamma_k)$, where $m_1,\dots,m_k \in \ZZ_{\geq 1}$, and the $\gamma_i$'s are pairwise disjoint closed geodesics on $X$ without self-intersections.
Our definition of a random (ordered) multi-geodesic is as follows.
Consider the set
\[
    s_{X,R,\gamma}
    \coloneqq
    \{ \alpha \in \Mod(X) \cdot \gamma : \ell_X(\alpha) \leq R \}
\]
of ordered multi-geodesics of the same topological type as that of $\gamma$ and of length at most $R$, where $\Mod(X)$ is the mapping class group of $X$, and $\ell_X(\alpha) \coloneqq m_1\ell_X(\alpha_1) + \cdots + m_k\ell_X(\alpha_k)$ is the total length of $\alpha$.
The set $s_{X,R,\gamma}$ is finite.
Consider the uniform probability measure on this set, and the random variable that associates a multi-geodesic to its normalized length vector
\[
    \hat{\ell}_{X,R,\gamma}
    \colon
    s_{X,R,\gamma}
    \to
    \varDelta^{k-1},
    \qquad
    (m_1 \alpha_1, \dots, m_k \alpha_k)
    \mapsto
    \frac{1}{\ell_{X}(\alpha)} \, (m_1\ell_X(\alpha_1), \dots, m_k\ell_X(\alpha_k))
\]
where $\varDelta^{k-1} \coloneqq \{ (x_1,\dots,x_k) \in \RR_{\geq 0}^k : x_1 + \cdots + x_k = 1\}$ is the standard simplex of dimension $k-1$.

Our goal is to find the limiting distribution of $\hat{\ell}_{X,R,\gamma}$ as $R\to\infty$.

The study of the length partition of random multi-geodesics is initialed by M.~Mirzakhani in \cite{Mir}, where she proves the following result.
\begin{theorem}[{\cite[Theorem~1.2]{Mir}}] \label{thm:Mir.1.2}
    If $\{ \gamma_1,\dots,\gamma_{3g-3} \}$ is a pants decomposition of $X$, then as $R\to\infty$, the random variable $\hat{\ell}_{X,R,\gamma}$ converges in law to the Dirichlet distribution of order $3g-3$ with parameters $(2,\dots,2)$, namely the probability distribution that has a density function $(6g-7)! \cdot x_1 \cdots x_{3g-3}$ with respect to the Lebesgue measure on the standard simplex $\varDelta^{3g-4} \coloneqq \{ (x_1,\dots,x_{3g-3}) \in \RR_{\geq 0}^{3g-3} : x_1 + \cdots + x_{3g-3} = 1 \}$ of dimension $3g-4$.
    In other words, for any open subset $U$ of $\varDelta^{3g-4}$,
    \[
        \lim_{R \to \infty} \PP(\hat{\ell}_{X,R,\gamma} \in U)
        =
        (6g-7)! \int_{U} x_1 \cdots x_{3g-3} \, \lambda(dx).
    \]
    where $\lambda$ is the Lebesgue measure on $\varDelta^{3g-4}$.
\end{theorem}
Our main result is the following generalization of the preceding theorem to any arbitrary topological type.
\begin{theorem} \label{thm:fixedTopo}
    Let $\gamma = (m_1\gamma_1,\dots,m_k\gamma_k)$ be an ordered multi-geodesic on $X$.
    As $R\to\infty$, the random variable $\hat{\ell}_{X,R,\gamma}$ converges in law to a random variable which admits a polynomial density with respect to the Lebesgue measure on $\varDelta^{k-1}$ given by, up to a normalizing constant,
    \[
        (x_1,\dots,x_k)
        \mapsto
        \bar{P}_{\gamma}(x_1/m_k,\dots,x_k/m_k)
    \]
    where $\bar{P}_{\gamma}$ is the top-degree (homogeneous) part of the graph polynomial $P_{\gamma}$ associated to $\gamma$ defined by \eqref{eq:Pgamma}.
\end{theorem}
\begin{remark}
    The function $\bar{P}_{\gamma}$ is a homogeneous polynomial of degree $6g-7$ whose coefficients can be expressed in terms of the psi-classes on the Deligne--Mumford compactification of the moduli space of smooth complex curves $\overline{\cM}_{g,n}$ (see Theorem~\ref{thm:V}).
    In particular, it depends only upon the topological type of $\gamma$, but not on the hyperbolic metric $X$.
\end{remark}

\subsection*{Motivations}
Theorem~\ref{thm:fixedTopo} is motivated by Theorem~\ref{thm:Mir.1.2} of Mirzakhani.
Another motivation originates from Theorem 1.25 in the section ``Statistical geometry of square-tiled surfaces'' of \cite{DGZZ21}.
Namely, the statistics of perimeters of maximal cylinders of a ``random'' square-tiled surface associated to a given multicurve $\gamma$ given by formula (1.38) from \cite{DGZZ21} coincide with statistics of hyperbolic lengths of different components of $\gamma$ given by Theorem \ref{thm:fixedTopo} above.
Though formula (1.38) from \cite{DGZZ21} can be interpreted as a certain \textit{average} of lengths statistics for individual hyperbolic surfaces $X$, it does not imply that such statistics do not change when $X$ changes.
The conjecture of non-varying of statistics of hyperbolic lengths for any fixed multicurve under arbitrary deformation of the hyperbolic metric, proved in Theorem~\ref{thm:fixedTopo}, was one of our principal motivations.

\subsection*{Idea of the proof} 
The structure of the proof is similar to that of \cite[Theorem~1.2]{Mir}.
The limiting distribution that we are after boils down to the asymptotics of multicurves counting under constraints, which can be transformed to a problem of approximating to the number of ``lattice points'' within a horoball in a covering space of the moduli space.
By considering tiling of the covering space by translates of a fundamental domain for the action of the mapping class group, it would not be unreasonable to expect that, this number might be proportional to the volume of the horoball divided by the volume of the moduli space, and finally this is not so far from the truth.
We proceed using techniques that Margulis introduced in his thesis \cite{Mar04}, and the equidistribution theorem for large horospheres initially established by Mirzakhani \cite{Mir07a}.
Similar methods were also applied in, e.g., \cite{EM93}.

Theorem~\ref{thm:fixedTopo} can be generalized to hyperbolic surfaces with cusps if Mirzakhani's work on the ergodicity of the earthquake flow can be generalized to such surfaces, which seems to be the case (see \cite{Mir}).

\subsection*{Remark}
While the author was finishing this paper, the paper \cite{AH21} by F.~Arana-Herrera appeared on the arXiv.
Both papers are devoted to a similar circle of problems and use a similar circle of ideas, though they were written in parallel and completely independently.
In particular, Arana-Herrera proves a much more general version of our Theorem~\ref{thm:equidisOnP1M} (\cite[Theorem 1.3]{AH21}), which is one of the key ingredients allowing to attack the counting problem and the length statistics.
We learned from \cite{AH21} that this kind of statistics was initially conjectured by S.~Wolpert.
Papers \cite{AH21} and \cite{AH22} established results closely related to Theorem~\ref{thm:fixedTopo}.

Proposition~\ref{prop:nuIsLebesgue} below is based on a theorem stated by M.~Mirzakhani but presented without a detailed proof.
The paper \cite{AH21} contains a detailed proof of an even stronger estimate which implies, in particular, the statements of this theorem; see Remark~\ref{rem:esti} below.

\subsection*{Acknowledgment}
The author is extremely grateful to Bram~Petri for all his patience and support, and for reading earlier drafts of this paper.
The author would also like to thank Anton~Zorich for formulating the problem, Francisco~Arana-Herrera for pointing out a gap in the first version of this paper, Grégoire~Sergeant-Perthuis and Jieao~Song for useful discussions, and the anonymous referee for a careful reading and numerous helpful suggestions.

\section{Background}
Throughout this paper, we use the symbol $\varSigma_{g}$ to denote a connected, closed, oriented, topological surface of genus $g \geq 2$, the symbol $\varSigma_{g,n}$ to denote a connected closed oriented topological surface of genus $g$ with $n$ boundary circles labeled by $\{ 1,\dots,n \}$ with $2g-2+n > 0$, and the letter $d$ to denote $3g-3+n$.

\subsection{Deformation spaces and mapping class group}
Let us consider the set of orientation-preserving homeomorphisms $\varphi \colon \varSigma_{g} \to X$ where $X$ is an oriented complete hyperbolic surface of genus $g$.
Two such homeomorphisms $\varphi_1 \colon \varSigma_{g} \to X_1$ and $\varphi_2 \colon \varSigma_{g} \to X_2$ are said to be \emph{equivalent} if $\varphi_2 \circ \varphi_1^{-1}$ is isotropic to an isometry.
The \emph{Teichmüller space}, denoted by $\cT(\varSigma_{g})$ or simply $\cT_{g}$, is the set of such equivalence classes.

Let us denote by $\Homeo^+(\varSigma_{g})$ the group of self-homeomorphisms of $\varSigma_{g}$ that preserve the orientation, and write $\Homeo_0(\varSigma_{g})$ for the subgroup of $\Homeo^+(\varSigma_{g})$ consisting of homeomorphisms isotropic to the identity. The \emph{mapping class group}, denoted by $\Mod(\varSigma_{g})$ or simply $\Mod_{g}$, is the quotient group $\Homeo^+(\varSigma_{g}) / \Homeo_0(\varSigma_{g})$.

The group $\Homeo^+(\varSigma_{g})$ acts (properly and discontinuously) from the right on $\cT_{g}$ by precomposition, and $\Homeo_0(\varSigma_{g})$ acts trivially. The \emph{moduli space}, denoted by $\cM(\varSigma_{g})$ or $\cM_{g}$, is the quotient $\cT_{g} / \Mod_{g}$.

The Teichmüller space $\cT_{g,n}(L_1,\dots,L_n)$ and moduli space $\cM_{g,n}(L_1,\dots,L_n)$ of oriented complete hyperbolic surfaces of genus $g$ with $n$ (labeled) totally geodesic boundary components of lengths $L_1, \dots, L_n \geq 0$ respectively can be defined in a similar manner.

\subsection{Curves}
In the introduction, the theorems are stated in terms of geodesics on a hyperbolic surface.
Nevertheless, it is often more convenient to work with (the free homotopy classes of) the topological curves, and they are actually equivalent for our purposes.
A curve in a topological space $X$ is (the image of) a continuous application $\bS^1 \to X$. In this paper, we are interested in curves up to free homotopy.
A closed curve is said to be \emph{simple} if it does not intersect itself.
A \emph{multicurve} is a finite multiset of disjoint simple curves, and a multicurve is \emph{ordered} (or \emph{labeled}) if its underlying set is labeled.
We will often write an ordered multicurve $\gamma$ as an ordered list $(m_1\gamma_1,\dots,m_k\gamma_k)$, and its unlabeled counterpart as a formal sum $\overline{\gamma} = m_1\gamma_1 + \cdots + m_k\gamma_k$ where $m_i \in \ZZ_{\geq 1}$, $1\leq i\leq k$.

The group $\Homeo^+(\varSigma_{g})$ acts on the set of closed curves on $\varSigma_{g}$ by postcomposition, and the action of the subgroup $\Homeo_0(\varSigma_{g})$ stabilizes sets of curves in the same free homotopy class. Thus the mapping class group acts on the set of free homotopy classes of closed curves on $\varSigma_{g}$.
We say that two closed curves $\alpha$ and $\beta$ have the same \emph{topological type} if they lie in the same mapping class group orbit.
The three following subgroups, associated to $\gamma$, of the mapping class group $\Mod_{g}$ will be useful later in the paper:
\begin{itemize}
    \item
        $\Stab(\overline{\gamma})$ which fixes the multicurve $\overline{\gamma} = m_1\gamma_1 + \cdots + m_k\gamma_k$ (but the $\gamma_i$'s can be permutated),

    \item
        $\Stab(\gamma)$ which fixes every $\gamma_i$ for all $1\leq i \leq k$,

    \item
        $\Stab^+(\gamma)$ which fixes every $\gamma_i$ and its orientation for all $1\leq i \leq k$.
\end{itemize}

Let $X \in \cT_{g}$.
If a closed curve $\alpha$ on $X$ is not homotopic to a point, then $\alpha$ is freely homotopic to a unique closed geodesic on $X$ with the minimum length over all curves in the free homotopy class of $\alpha$, and we write $\ell_X(\alpha)$ for the length of this geodesic.

The notions of topological type and length extend naturally to multicurves.

\subsection{Fenchel--Nielsen coordinates}
A \emph{pair of pants} is a surface that is homeomorphic to a sphere with three holes.
A \emph{pants decomposition} of $\varSigma_{g,n}$ is a set of disjoint simple closed curves $\{ \alpha_1,\dots,\alpha_{3g-3+n} \}$ on $\varSigma_{g,n}$ such that $\varSigma_{g,n} \smallsetminus \{ \alpha_1, \dots, \alpha_{3g-3+n} \}$ is a disjoint union of pairs of pants. 

Fix an ordered pants decomposition $(\alpha_1, \dots, \alpha_{3g-3+n})$ of $\varSigma_{g,n}$.
Given $X \in \cT_{g,n}(L_1,\dots,L_n)$ (or $X \in \cT_{g}$), we can associate for each $\alpha_i$ two parameters: the length $\ell_{\alpha_i}(X) \in \RR_{>0}$, and the \emph{twist parameter} $\tau_{\alpha_i}(X) \in \RR$ (corresponding to how much one turns before gluing two pairs of pants along $\alpha_i$; see \cite[Section 1.7]{Bus92} for a precise definition). These $6g-6+2n$ parameters are called \emph{Fenchel--Nielsen coordinates}. The application $\cT_{g,n}(L_1,\dots,L_n) \to (\RR_{>0} \times \RR)^{3g-3+n}$ given by $X \mapsto (\ell_{\alpha_i}(X), \tau_{\alpha_i}(X))_{i=1}^{k}$ is a bijection (see \cite[Chapter 6]{Bus92}).

\subsection{Weil--Petersson volumes}
The following theorem is often referred to as Wolpert’s magical formula.
\begin{theorem}[{\cite{Wol83}}]
    Given a pants decomposition $\{ \alpha_1,\dots,\alpha_{3g-3+n} \}$, the formula
    \[
        \sum_{i=1}^{3g-3+n} d\ell_{\alpha_i} \wedge d\tau_{\alpha_i}
    \]
    defines a symplectic form which has the same expression in any other Fenchel--Nielsen coordinates.
    In particular, it is invariant under the action of the mapping class group.
\end{theorem}
The symplectic form thus defined is the so-called \emph{Weil--Petersson symplectic form}, and we shall denote it by $\omega$.
See \cite{Wol83} for a more intrinsic definition.

Every symplectic form defines a volume form. The \emph{Weil--Petersson volume} of the moduli space $\cM_{g,n}(L_1,\dots,L_n)$ is defined by
\[
    V_{g,n}(L_1,\dots,L_n)
    \coloneqq
    \int_{\cM_{g,n}(L_1,\dots,L_n)} \frac{\omega^{\mkern1mu \wedge (3g-3+n)}}{(3g-3+n)!}.
\]
The following fundamental result is due to Mirzakhani.
\begin{theorem}[{\cite{Mir07c}}] \label{thm:V}
    The Weil--Petersson volume $V_{g,n}(L_1,\dots,L_n)$ is a symmetric polynomial in $L_1^2, \dots, L_n^2$ of degree $3g-3+2n$. More precisely,
    \[
        V_{g,n}(L_1,\dots,L_n)
        =
        \sum_{\substack{(d_0,d_1,\dots,d_n) \in \ZZ_{\geq 0} \\ d_0 + d_1 + \cdots + d_n = 3g-3+n}} \frac{(2\pi^2)^{d_0}}{2^{d_1 + \cdots + d_n} d_0! d_1! \cdots d_n!} \left( \int_{\overline{\cM}_{g,n}} \kappa_1^{d_0} \psi_1^{d_1} \cdots \psi_n^{d_n} \right) L_1^{2d_1} \cdots L_n^{2d_n}
    \]
    where $\overline{\cM}_{g,n}$ is the Deligne--Mumford compactification, $\psi_i \in H^2(\overline{\cM}_{g,n},\QQ)$ is the $i$-th psi-class, and $\kappa_1 = [\omega]/2\pi^2 \in H^2(\overline{\cM}_{g,n},\QQ)$ is the first Mumford class.
\end{theorem}

\subsection{Earthquakes}
Multicurves can be regarded as ``lattice points'' in the space of \emph{measured laminations} $\cML_{g}$.
We will only need the following properties of this space.
See, e.g., \cite[Chapter~11]{Kap01} for more details.
\begin{enumerate}
    \item[1.]
        The space $\cML_{g}$ is a $(6g-6)$-dimensional real manifold equipped with a natural piece-wise integral linear structure, i.e., $\cML_{g}$ has an natural atlas whose transition functions are piece-wise in $\GL(6g-6,\ZZ)$.

    \item[2.]
        The integral points in the coordinate charts of $\cML_{g}$, denoted by $\cML_{g}(\ZZ)$, are in natural bijection with the (free homotopy classes of) integral multicurves on $\varSigma_{g}$.

    \item[3.]
        The action of the mapping class group on the set of multicurves extends to $\cML_{g}$.

    \item[4.]
        The group $(\RR_{>0},\times)$ acts on $\cML_{g}$, and for any multicurve $\overline{\gamma} = m_1\gamma_1 + \cdots + m_k\gamma_k$, and any $r\in\RR_{>0}$, $r \cdot (m_1\gamma_1 + \cdots + m_k\gamma_k) = r\, m_1\gamma_1 + \cdots + r\, m_k\gamma_k$. We denote the quotient by $\PP(\cML_{g})$.

    \item[5.]
        Given $X \in \cT_{g}$, the length function $\ell_X$ defined on the set of multicurves extends to $\cML_{g}$.
        Moreover, for any $\lambda \in \cML_{g}$, we have $\ell_{X\cdot h}(h^{-1} \cdot \lambda) = \ell_X(\lambda)$ for any $h \in \Mod_{g}$, and $\ell_X(r \cdot \lambda) = r \cdot \ell_X(\lambda)$ for any $r\in \RR_{>0}$.

    \item[6.]
        The twist flow $\tw_{\overline{\gamma}}^{\mkern2mu t}$ about a multicurve $\overline{\gamma}$ can be extended to any measured lamination $\lambda \in \cML_{g}$, and we have $(\tw_{\lambda}^{\mkern2mu t}(X)) \cdot h = \tw_{h^{-1} \lambda}^{\mkern2mu t}(X h)$ and $\tw_{r \lambda}^{\mkern2mu t}(X) = \tw_{\lambda}^{\mkern2murt}$ for all $t\in \RR$, $h \in \Mod_{g}$, and $r\in \RR_{>0}$.

    \item[7.]
        The space $\cML_{g}$ carries a natural mapping class group invariant measure $\mu_{\Th}$ defined by asymptotic counting of integral points, called the \emph{Thurston measure}. The Thurston measure is a Lebesgue measure in the coordinate charts of $\cML_g$, and for any open subset $U \subset \cML_{g}$, we have $\mu_{\Th}(t\cdot U) = t^{6g-6} \mu_{\Th}(U)$ for any $t\in \RR_{>0}$.
\end{enumerate}

Let $\cP\cT_{g} \coloneqq \cT_{g} \times \cML_{g}$ be the bundle of measured laminations over the Teichmüller space, and let $\cP^1\cT_{g} \coloneqq \{ (X,\lambda) \in \cP\cT_{g} : \ell_X(\lambda) = 1\}$ be the unit sphere bundle of $\cP\cT_{g}$ with respect to the length function.

The mapping class group acts on $\cP\cT_{g}$ from the right via $(X,\lambda) \cdot h \coloneqq (X \cdot h, h^{-1} \cdot \lambda)$.
This action is well-defined on $\cP^1\cT_{g}$ since it preserves the length function $\ell(X,\lambda) \coloneqq \ell_X(\lambda)$. Write $\cP\cM_{g} \coloneqq \cP\cT_{g} / \Mod_{g}$ and $\cP^1\cM_{g} \coloneqq \cP^1\cT_{g} / \Mod_{g}$.

The \emph{earthquake flow} $\tw^{\mkern1mu t}$ on $\cP\cT_{n}$ is defined by
\[
    \tw^{\mkern2mu t}(X,\lambda)
    \coloneqq
    (\tw^{\mkern2mu t}_\lambda(X),\lambda).
\]
The earthquake flow commutes with the action of the mapping class group, and therefore descends to $\cP\cM_{g}$, and to $\cP^1\cM_{g}$ (since the earthquake preserves the length function).

The Thurston measure on $\cML_{g}$ induces a measure on $\{ \lambda \in \cML_{g} : \ell_X(\lambda) = 1 \}$ in the following way: let $U\subset \{ \lambda \in \cML_{g} : \ell_X(\lambda) = 1 \}$ be an open subset. The \emph{Thurston measure} of $U$ is defined to be 
\[
    \mu_\Th\{ s \cdot \lambda \in \cML_{g} : \lambda \in U,\ s\in [0,1]\}.
\]
The measure $\nu_{g}$ on $\cP^1\cT_{g}$ defined by
\[
    \nu_{g}(U)
    \coloneqq
    \int_{\cT_{g}} \mu_{\Th}\{ s \lambda \in \cML_{g} : (X,\lambda) \in U,\ s\in [0,1] \} \, dX
\]
for any open subset $U \subset \cP^1\cT_{g}$, is invariant both  under the earthquake flow (since $\mu_\WP$ is) and under the action of the mapping class group (since $\mu_\Th$ and $\mu_\WP$ are), and hence descends to a measure on $\cP^1\cM_{g}$ that (by abuse of notation) we shall also denote by $\nu_g$. The total mass of $\nu_g$
\[
    b_g
    =
    \int_{\cM_g} B(X) \, dX
\]
where $B(X) \coloneqq \mu_\Th \{ \lambda \in \cML_{g} : \ell_X(\lambda) \leq 1 \}$, is finite \cite[Theorem~3.3]{Mir08a}.

The following result is fundamental.
\begin{theorem}[\cite{Mir08a}]
    The earthquake flow on $\cP^1\cM_{g}$ is ergodic with respect to $\nu_{g}$.
\end{theorem}
We recommend \cite{Wri} for an expository survey on this topic.
\subsection{Thurston distance}
Let $X_1,X_2\in\cT_{g}$. Set
\[
    d(X,Y)
    \coloneqq
    \sup_{\lambda \in  \cML_{g}} \log\frac{\ell_{X_1}(\lambda)}{\ell_{X_2}(\lambda)}.
\]
The \emph{Thurston distance} between $X_1$ and $X_2$ is defined by
\[
    d_\Th(X_1,X_2)
    \coloneqq
    \max\{ d(X_1,X_2), d(X_1,X_2) \}.
\]
The \emph{Thurston distance ball} centered at $X \in \cT_{g}$ of radius $\epsilon$ is defined to be
\[
    \BB_X(\epsilon)
    \coloneqq
    \{ Y \in \cT_{g} : d_\Th(X,Y) \leq \epsilon/2\}.
\]
The reason for this choice of radius is that, for small $\epsilon$, e.g.\ $0<\epsilon < 1$,
\[
    e^\epsilon < 1 + 2x,
    \qquad
    e^{-\epsilon} > 1 - 2x.
\]
We have therefore, for any $\lambda \in \cML_{g}$ and any $Y \in \BB_X(\epsilon)$,
\begin{equation} \label{eq:1+-e}
    (1-\epsilon) \cdot \ell_X(\lambda)
    \leq
    \ell_Y(\lambda)
    \leq
    (1+\epsilon) \cdot \ell_X(\lambda)
\end{equation}
Thurston distance balls are well-defined on $\cM_{g}^\gamma$, and on $\cM_{g}$, since the Thurston distance is $\Mod_{g}$-invariant.

\subsection{Stable graphs}
Given a multicurve $m_1\gamma_1 + \cdots + m_k\gamma_k$, one can associate with it a \emph{stable graph} in the following way.
Cut the surface along $\gamma_1,\dots,\gamma_k$.
To each connected component $S$ of $\varSigma_{g} \smallsetminus \{ \gamma_1,\dots,\gamma_k \}$, we associate a vertex, and we decorate this vertex with the genus of $S$.
For each component $\gamma_i$ of $\gamma$, we draw an edge that connects the two vertices (which could be the same) corresponding to the two connected components of $\varSigma_{g} \smallsetminus \{ \gamma_1,\dots,\gamma_k \}$ bounded by $\gamma_i$.
See Figure~\ref{fig:1} for an example.
Note that the resulting graph does not depend on $m_1,\dots,m_k$.
More formally, a stable graph consists of the data
\[
    \varGamma
    =
    (V,\ E,\ H,\ g \colon V\to\ZZ_{\geq 0},\ \iota \colon H \to H)
\]
satisfying the following properties:
\begin{enumerate}
    \item[1.]
        The pair $(V,E)$ defines a connected graph, with vertex set $V$ and edge set $E$. The set $H$ is the set of \emph{half-edges}.

    \item[2.]
        The map $v$ assigns each half-edge to its adjacent vertex.

    \item[3.]
        The map $\iota$ is an involution, such that the $2$-cycles of $\iota$ are in bijection with $E$, and the fixed points of $\iota$ are in bijection with $L$.

    \item[4.]
        The \emph{genus function} $g$ assigns each vertex $x$ to its genus (the genus of the surface corresponding to $x$), such that the \emph{stability condition}
        \[
            2g(x) - 2 + n(x) > 0
        \]
        is satisfied, where $n(x)$ denotes the number of edges and legs adjacent to $x$.
\end{enumerate}

\subsection{Graph polynomials}
Given a stable graph associated to the multicurve $\gamma = m_1\gamma_1 + \cdots + m_k\gamma_k$, we associate to each edge $e$ a variable $x_e$, and define the associated \emph{graph polynomial} by the formula
\begin{equation} \label{eq:Pgamma}
    P_\gamma(x_e : e \in E)
    =
    \prod_{e} x_e \cdot \prod_{v} V_{g(v),n(v)}(x_{e(h)} : h\in H,\ v(h) = v).
\end{equation}
where $e$ runs through the edge set $E$, $v$ runs through the vertex set $V$, $V_{g(v),n(v)}$ is the Weil--Petersson volume polynomial of $\cM_{g(v),n(v)}$ (see Theorem~\ref{thm:V}), $e(h)$ is the edge that contains the half-edge $h$, and $v(h)$ denotes the vertex incident to $h$.
Note that $P_{\gamma}$ is of degree $2d-k$.
Finally, we write $\bar{P}_\gamma$ for the top-degree homogeneous part of $P_{\gamma}$, and $\bar{V}_{g,n}$ for that of $V_{g,n}$.

\begin{example}
    If $\{ \gamma_1,\dots,\gamma_{3g-3} \}$ is a pants decomposition, then $V_{g(v),n(v)} = 1$ for all $v\in V$, and $P_\gamma(x_1,\dots,x_{3g-3}) = \bar{P}_\gamma(x_1,\dots,x_{3g-3}) = x_1\cdots x_{3g-3}$.
\end{example}

\begin{example} \label{ex:1}
    Let $(\gamma_1, \gamma_2, \gamma_3)$ be an ordered multicurve on $\varSigma_{3}$ as in Firgure~\ref{fig:1}.
    The Weil--Petersson volume polynomial $V_{1,3}(x_1,x_2,x_3)$ is equal to (see \cite{Do15})
    \[
        \left( \frac{\msp_{(3)}}{1152} + \frac{\msp_{(2,1)}}{192} + \frac{\msp_{(1,1,1)}}{96} + \frac{\pi^2 \msp_{(2)}}{24} + \frac{\pi^2 \msp_{(1,1)}}{8} + \frac{13\pi^4 \msp_{(1)}}{24}\right) (x_1^2,x_2^2,x_3^2) + \frac{14 \pi^6}{9}
    \]
    where $\msp$ stands for the monomial symmetric polynomial. For example,
    \[
        \msp_{(2,1)}(x_1,x_2,x_3)
        =
        x_1^2 x_2 + x_1 x_3^2 + x_1 x_2^2 + x_2^2 x_3 + x_1 x_3^2 + x_2 x_3^2,
        \quad
        \msp_{(1)}(x_1,x_2,x_3)
        =
        x_1 + x_2 + x_3.
    \]
    So the top-degree part of $V_{1,3}$ is
    \[
        \bar{V}_{1,3}(x_1,x_1,x_2)
        =
        \frac{2x_1^6 + x_2^6}{1152} + \frac{2x_1^6 + 2x_1^4x_2^2 + 2x_1^2x_2^4}{192} + \frac{x_1^4 x_2^2}{96},
    \]
    and therefore,
    \[
        \bar{P}_{\gamma}(x_1,x_2,x_3)
        =
        \frac{2x_1^6 x_2 x_3 + x_1 x_2^6 x_3}{1152} + \frac{x_1^7x_2x_3 + x_1^5x_2^3x_3 + x_1^3 x_2^5 x_3 + x_1^5 x_2^3 x_3}{96}.
    \]
    \begin{figure}[h] 
        \centering
        \begin{overpic}[width = 0.6\textwidth]{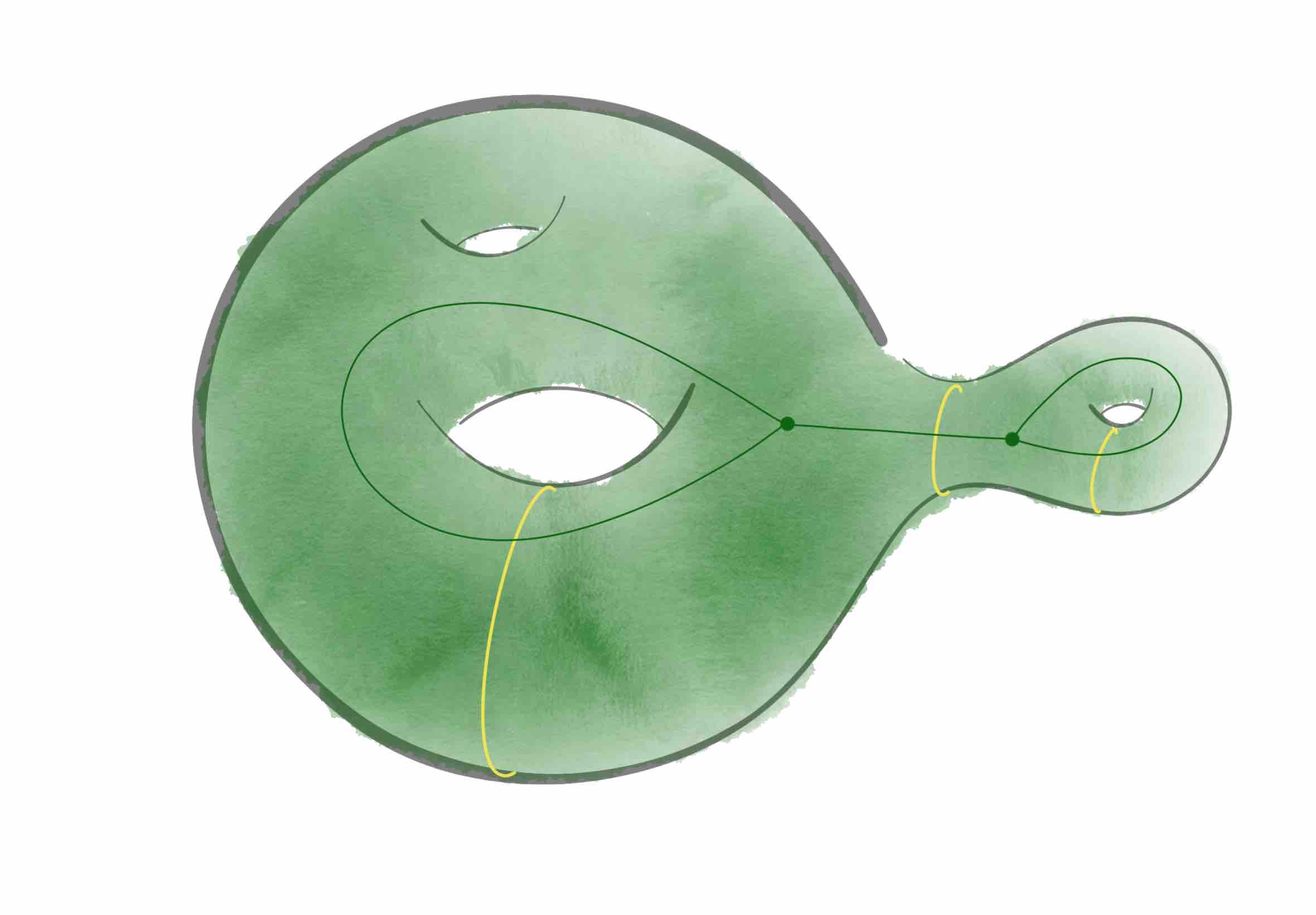}
            \put(36,7){$\gamma_1$}
            \put(70,28){$\gamma_2$}
            \put(82,27){$\gamma_3$}
        \end{overpic}
        \caption{Example~\ref{ex:1}} \label{fig:1}
    \end{figure}
\end{example}

\section{Mirzakhani's covering spaces} \label{sec:Mggamma}
Let $\gamma = (m_1\gamma_1, \dots, m_k\gamma_k)$ be an ordered multicurve on $\varSigma_{g}$.
Recall that $\Stab(\gamma)$ denotes the subgroup of $\Mod_g$ that fixes each $\gamma_i$, $1 \leq i \leq k$.
The quotient space 
\[
    \cM_{g}^{\gamma} \coloneqq \cT_{g} / \Stab(\gamma)
\]
introduced by Mirzakhani in her thesis plays an important role in this paper.

Write $\pi^{\gamma} \colon \cT_{g} \to \cM_{g}^{\gamma}$ and $\pi_{\gamma} \colon \cM_{g}^\gamma \to \cM_{g}$ for the two natural projections (``raising and lowering the index'').
Let us consider the product space $P_\gamma \coloneqq \cT_{g} \times \Mod_{g} \cdot \, (\gamma_1,\dots,\gamma_k)$ of the Teichmüller space and the mapping class group orbit of $(\gamma_1, \dots, \gamma_k)$.
The mapping class group $\Mod_{g}$ acts on $P$ (from the right) via $(X ; \alpha_1, \dots, \alpha_k) \cdot h = (X \cdot h ; h^{-1} \alpha_1, \dots, h^{-1} \alpha_k)$.

\begin{lemma}
    The quotient $P_\gamma / \Mod_{g}$ is isomorphic to $\cM_{g}^{\gamma}$ as symplectic orbifolds.
\end{lemma}
\begin{proof}
    Consider the map $P_{\gamma} \to \cM_{g}^\gamma$ defined by $(X, h \gamma) \mapsto \pi^{\gamma}(X h)$.
    This map is surjective, and descends to the quotient $P_\gamma/\Mod_{g}$.
    The resulting map $P_\gamma/\Mod_g \to \cM_{g}^{\gamma}$ is a local isomorphism of symplectic orbifolds.
    All that remains now is to show that the map $P_\gamma/\Mod_{g} \to \cM_{g}^\gamma$ is injective.
    Let $(X_1,h_1 \gamma)$, $(X_2, h_2 \gamma) \in P_\gamma$ such that $\pi^\gamma(X_1 h_1) = \pi^\gamma(X_2 h_2)$.
    By definition, there exists $s\in\Stab(\gamma)$ such that $X_1 h_1s = X_2 h_2$.
    Therefore
    \[
        (X_2, h_2 \gamma)
        =
        (X_1 h_1sh_2^{-1}, h_2 \gamma)
        \sim
        (X_1 h_1s, \gamma)
        \sim
        (X_1 h_1, \gamma),
    \]
    which proves the injectivity.
\end{proof}
\begin{remark} \label{rem:lXgamma}
    Let $\alpha$ be a simple closed curve on $\varSigma_{g}$.
    In general, $\ell_X(\alpha)$ is not well-defined for $X \in \cM_{g}^\gamma$.
    However, it is if $\alpha = \gamma_i$ for some $i$.
\end{remark}
The next lemma is a simple fact, but for our purposes it will be very important: it transforms the multicurves counting that we are after to a ``lattice points'' counting problem on $\cM_{g}^{\gamma}$.
\begin{lemma} \label{lem:bij}
    Let $\gamma = (m_1\gamma_1,\dots,m_k\gamma_k)$ be an ordered multicurve, $X\in\cT_{g}$, $R \in \RR_{>0}$, and $A \subset \varDelta^{k-1}$ be an open subset.
    The set
    \[
        \{ h \gamma : h \in \Mod_{g},\ \ell_X(h \gamma) \leq R,\ \hat{\ell}_X(h \gamma) \in A \}
    \]
    and the set
    \[
        \{ [(X,h \gamma)] \in \cM_{g}^\gamma : h \in \Mod_{g},\ \ell_X(h\gamma) \leq R,\ \hat{\ell}_X(h \gamma) \in A\},
    \]
    where by $[(X, h \gamma)]$ we mean the image of $(X,h\gamma)$ under $P \to P/\Mod_{g}$, are in bijection given by $h \gamma \mapsto [(X, h \gamma)]$,
\end{lemma}
\begin{proof}
    The given map is obviously surjective.
    Suppose that $[(X, h_1 \gamma)] = [(X, h_2 \gamma)]$, then $h_1^{-1}h_2 \in \Stab(\gamma)$, and therefore $h_1\gamma = h_2\gamma$. The injectivity follows. 
\end{proof}

Next, let us review another covering space of $\cM_{g}$ that Mirzakhani introduced.
By considering the Fenchel--Nielsen coordinates associated to a pants decomposition that contains $\gamma_1,\dots,\gamma_k$, the Teichmüller space $\cT_{g}$ can be written as
\begin{equation} \label{eq:Tgn}
    \{ (\ell_e, \tau_e, X_v) : e \in E,\ v \in V,\ \ell_e \in \RR_{>0},\ \tau_e \in \RR,\ X_v\in \cT_{g(v),n(v)}(\ell_{e(h)} : h\in H,\ v(h) = v)\}
\end{equation}
where $V$ (resp.\ $E$; $H$) is the vertex (resp.\ edge; half-edge) set of the stable graph associated to $\gamma$.
The group
\[
    G_{\gamma}
    \coloneqq
    \prod_{e} \ZZ \times \prod_{v} \Mod_{g(v),n(v)}
\]
acts naturally on $\cT_{g}$ written in the form \eqref{eq:Tgn} (each copy of $\ZZ$ acts as the Dehn twist about a $\gamma_i$), and $G_{\gamma}$ can be identified with $\Stab^+(\gamma)$. The quotient $C_\gamma \coloneqq \cT_{g} / G_\gamma$ is of the form
\[
    \{ (\ell_e, \tau_e, X_v) : e \in E,\ v \in V,\ \ell_e \in \RR_{>0},\ \tau_e \in \RR/\ell_e \,\ZZ,\ X_v\in \cM_{g(v),n(v)}(\ell_{e(h)} : h \in H,\ v(h) = v)\}.
\]
Since $G_\gamma \simeq \Stab^+(\gamma)$ is a subgroup of $\Stab(\gamma)$, $\cT_{g} \to \cM_{g}^{\gamma}$ factors through a (ramified) covering map $C_{\gamma} \to \cM_{g}^{\gamma}$. The degree of this covering map is
\[
    \kappa_\gamma
    =
    2^{\mkern1mu M(\gamma)} \cdot [\Stab(\gamma) : \langle \Stab^+(\gamma), \Stab_0(\gamma) \rangle]
\]
where $M(\gamma)$ is the number of $i$ such that $\gamma_i$ bounds a surface homeomorphic to $\varSigma_{1,1}$, and $\langle \Stab^+(\gamma), \Stab_0(\gamma) \rangle$ stands for the subgroup of $\Stab(\gamma)$ generated by $\Stab^+(\gamma)$ and the kernel $\Stab_0(\gamma)$ of the action of $\Stab(\gamma)$ on $\cT_{g}$.
Note that $\Stab_0(\gamma)$ is trivial when $g \geq 3$, and is isomorphic to $\ZZ / 2\ZZ$ if $g=2$ (generated by the hyperelliptic involution which fixes the free homotopy class of every simple closed curve on $\varSigma_{2}$).
For more details, see the footnote on p.~369--370 of \cite{Wri20}.

Integrating functions over $C_\gamma$ (and $\cM_{g}^{\gamma}$) is far less delicate than integrating function over $\cM_{g}$.
Starting from this observation Mirzakhani was able to calculate the integrals of an important class of functions defined on $\cM_{g}$, which she called ``geometric functions''.
\begin{theorem}[Mirzakhani's integration formula] \label{thm:MirIntFor}
    Let $\gamma = (\gamma_1,\dots,\gamma_k)$ be an ordered multicurve, and $f\colon \RR_{>0}^{k} \to \RR$ be a measurable function. Let $X \in \cM_{g}$, and choose an $\widetilde{X} \in \pi^{-1}(X) \in \cT_{g}$. We define $f_\gamma \colon \cM_{g} \to \RR$ by the formula
    \[
        f_\gamma(X)
        \coloneqq
        \sum_{(\alpha_1,\dots,\alpha_k) \in \Mod_{g} \cdot (\gamma_1,\dots,\gamma_k)} f(\ell_{\widetilde{X}}(\alpha_1), \dots, \ell_{\widetilde{X}}(\alpha_k))
    \]
    Note that $f_\gamma(X)$ does not depend on the choice of $\widetilde{X}$.
    We have
    \begin{align*}
        \int_{\cM_{g}} f_\gamma(X) \, dX
        & =
        \int_{\cM_{g}^\gamma} f(\ell_{X}(\gamma_1), \dots, \ell_{X}(\gamma_k)) \, dX \\
        & =
        \kappa_\gamma \int_{\RR_{>0}^k} f(x_1,\dots,x_k) \cdot P_\gamma(x_1,\dots,x_k) \, dx_1\cdots dx_k.
    \end{align*}
\end{theorem}

\section{Horospheres}
Let $\gamma = (m_1\gamma_1,\dots,m_k\gamma_k)$ be an ordered multicurve, $\overline{\gamma} = m_1\gamma_1 + \cdots + m_k\gamma_k$ be its unlabeled counterpart, and $A$ be an open subset of the standard simplex $\varDelta^{k-1} \coloneqq \{ (x_1,\dots, x_k) \in \RR_{\geq 0}^k : x_1 + \cdots + x_k = 1 \}$ of dimension $k-1$.

The \emph{horosphere} of radius $R$ associated to $\gamma$ and $A$ on $\cT_{g}$ is defined by
\[
    \tilde{\cS}_{R, \gamma}^{A}
    \coloneqq
    \{ X \in \cT_{g} : \ell_X(\gamma) = R,\ \hat{\ell}_X(\gamma) \in A \}.
\]
Similar notions can be defined on $\cM_{g}^{\gamma}$ and on $\cM_{g}$ by
\[
    \cS_{R, \gamma}^{A}
    \coloneqq
    \pi^\gamma(\tilde{\cS}_{R, \gamma}^{A}) \subset \cM_{g}^\gamma,
    \qquad
    \bar{\cS}_{R, \gamma}^{A}
    \coloneqq
    \pi(\tilde{\cS}_{R, \gamma}^{A}) \subset \cM_{g}
\]
where $\pi^{\gamma} \colon \cT_{g} \to \cM_{g}^{\gamma}$ and $\pi \colon \cT_{g} \to \cM_{g}$ are the natural projections. 
\begin{remark} \label{rem:l-1}
    The function $X \mapsto (m_i\ell_X(\gamma_i))_{i=1}^{k}$ is well-defined for $X \in \cM_{g}^{\gamma}$.
    The horosphere $\cS_{R, \gamma}^{A} \subset \cM_{g}^{\gamma}$ can be written as the pre-image of $R \cdot A$ under this function, where $R \cdot A$ is defined to be $\{ (x_1,\dots,x_k) \in \RR_{>0}^k : (x_1,\dots,x_k)/R \in A \}$.
\end{remark}

\subsection{Horospherical measures}
We can choose $d-k$ simple closed curves $\alpha_{k+1},\dots,\alpha_{d}$ such that $\{ \gamma_1,\dots,\gamma_k,\alpha_{k+1},\dots,\alpha_d \}$ is a pants decomposition.
In the associated Fenchel--Nielsen coordinates, the horosphere $\tilde{\cS}_{R, \gamma}^{A}$ is an open subset of a simplex.
Let $\mu_{\varDelta}$ denote the Weil--Petersson (Lebesgue) measure on this simplex.
The \emph{horospherical measure} $\mu_{R,\gamma}^{A}$, of an open subset $U\subset\cT_{g}$ is defined to be
\[
    \mu_{R,\gamma}^{A}(U)
    \coloneqq
    \mu_{\varDelta}(U \cap \tilde{\cS}_{R,\gamma}^{A}).
\]
The horospherical measure $\mu_{R,\gamma}^{A}$ is invariant under the action of the mapping class group, and hence descends to a measure on $\cM_{g}^\gamma$ and a measure on $\cM_{g}$; by abuse of notation we shall denote both by $\mu_{R,\gamma}^{A}$.
Note that $\cM_{g}^{\gamma} \to \cM_{g}$ is a (ramified) covering map of infinite degree.
However, its restriction on $\cS_{R,\gamma}^{A}$ is of finite degree.
Thus $\mu_{R,\gamma}^{A}$ on $\cM_{g}$ is the push-forward measure of $\mu_{R,\gamma}^{A}$ by $\cM_{g}^{\gamma} \to \cM_{g}$. So for any open subset $U$ of $\cM_{g}$,
\[
    \mu_{R,\gamma}^{A}(\pi_{\gamma}^{-1}(U))
    =
    [\Stab(\overline{\gamma}) : \Stab(\gamma)] \cdot \mu_{R,\gamma}^{A}(U).
\]
In particular, the total masses of $\mu_{R,\gamma}^{A}$ on $\cM_{g}^{\gamma}$ and on $\cM_{g}$ differ only by a multiplicative constant depending only on $\gamma$.
\subsection{Total mass}
The horospherical measure on $\cT_{g}$ has infinite total mass.
Nevertheless, its total mass is finite on $\cM_{g}^{\gamma}$ and $\cM_{g}$.
\begin{proposition}
    The total mass of $\mu_{R,\gamma}^{A}$ on $\cM_{g}$ is
    \[
        M_{R,\gamma}^{A}
        =
        \frac{\kappa_\gamma}{[\Stab(\overline{\gamma}) : \Stab(\gamma)]}\frac{1}{m_1\cdots m_k} \int_{R\cdot A} P_{\gamma}(x_1/m_k,\dots,x_k/m_k)\, \lambda(dx)
    \]
    where $R\cdot A \coloneqq \{ (x_1,\dots,x_k) \in \RR_{\geq 0}^k : (x_1,\dots,x_k)/R \in A \}$, $\lambda$ is the Lebesgue measure on $\{ (x_1,\dots,x_k) \in \RR_{\geq 0} : x_1 + \cdots + x_k = R\}$, and $P_\gamma$ is defined by the formula \eqref{eq:Pgamma}.
\end{proposition}
\begin{proof}
    In the light of Remark~\ref{rem:l-1}, by taking $f$ in Theorem~\ref{thm:MirIntFor} to be the indicator function
    \[
        \mathds{1}\left\{ (x_1,\dots,x_k) \in \RR_{>0}^k :
            \begin{array}{lr}
                R \leq m_1x_1 + \cdots + m_kx_k \leq R + \epsilon, \\
                (m_1x_1,\dots,m_kx_k)/(m_1x_1 + \cdots +m_kx_k) \in A
        \end{array} \right \}
    \]
    we obtain that $\mu_{R,\gamma}^{A}(\cM_{g}) \cdot [\Stab(\overline{\gamma}) : \Stab(\gamma)]$ is equal to
    \[
        \mu_{R,\gamma}^{A}(\cM_{g}^{\gamma})
        =
        \lim_{\epsilon \to 0} \frac{1}{\epsilon} \int_{\cM_{g}} f_{\gamma}(X) \, dX
        =
        \frac{\kappa_\gamma}{m_1\cdots m_k}\int_{R \cdot A} P_{\gamma}(x_1/m_1,\dots,x_k/m_k)\, \lambda(dx),
    \]
    the result desired.
\end{proof}
\begin{corollary} \label{cor:M}
    The total mass $M_{R,\gamma}^{A}$ (of $\mu_{R,\gamma}^{A}$ on $\cM_{g}$) is a polynomial in $R$ of degree $2d-1 = 6g-7$. Write $C_{\gamma}^{A}$ for its leading coefficient. We have
    \begin{equation} \label{eq:M.=CR2d-1}
        M_{R,\gamma}^{A}
        \sim
        C_{\gamma}^{A} \cdot R^{2d-1}
    \end{equation}
    as $R\to\infty$, and $C_{\gamma}^{A}$ can be calculated by
    \[
        C_{\gamma}^{A}
        =
        \frac{\kappa_\gamma}{[\Stab(\overline{\gamma}) : \Stab(\gamma)]} \frac{1}{m_1 \cdots m_k} \int_{A} \bar{P}_{\gamma}(x_1/m_1,\dots,x_k/m_k) \, \lambda(dx),
    \]
    where $\lambda$ is the Lebesgue measure on $\varDelta^{k-1}$ and $\bar{P}_\gamma$ is the top-degree homogeneous part of the graph polynomial $P_{\gamma}$ defined by \eqref{eq:Pgamma}.
\end{corollary}
\begin{remark}
    It results from Theorem~\ref{thm:V} that the polynomial $\bar{P}_\gamma$ can be expressed in terms of intersections numbers of $\psi$-classes on the Deligne--Mumford compactification $\overline{\cM}_{g,n}$.
\end{remark}
\subsection{Horospherical measures on the unit sphere bundle}
Define
\[
    \cP\tilde{\cS}_{R,\gamma}^{A}
    \coloneqq
    \{ (X, \gamma/R) \in \cP^1\cT_{g} : \hat{\ell}_X(\gamma) \in A \}.
\]
Note that $\cP\tilde{\cS}_{R,\gamma}^{A}$ projects via $\cP^1\cT_{g} \to \cT_{g}$ to $\tilde{\cS}_{R,\gamma}^{A}$, and is invariant under the earthquake flow.
Let $\nu_\varDelta$ denote the Lebesgue measure on $\cP\tilde{\cS}_{R,\gamma}^{A}$. The \emph{horospherical measure} $\nu_{R,\gamma}^{\mkern1mu A}$ on $\cP\cT_{g}$ is defined by the formula
\[
    \nu_{R,\gamma}^{\mkern1mu A}(U)
    \coloneqq
    \nu_{\varDelta}(U \cap \cP\tilde{\cS}_{R,\gamma}^{A})
\]
where $U$ is any open subset of $\cP^1\cT_{g}$. The measure $\nu_{R,\gamma}^{A}$ is $\Mod_{g}$-invariant, and therefore descends to a measure on $\cP^1\cM_{g}$ which by abuse of notation we shall also denote by $\nu_{R,\gamma}^{A}$. Note that $\mu_{R,\gamma}^{A}$ is the push-forward of $\nu_{R,\gamma}^{A}$ via $\cP^1\cM_{g} \to \cM_{g}$.

\subsection*{Notation}
To simplify the notation, let us fix $X \in \cM_{g}$, a multicurve $\gamma = (m_1\gamma_1, \dots, m_k\gamma_k)$ on $\varSigma_{g}$, and an open subset $A$ of the standard simplex $\varDelta^{k-1} \coloneqq \{ (x_1,\dots,x_k) \in \RR_{\geq 0} : x_1 + \cdots + x_k = 1 \}$.
From now on we shall write $\mu_R$ for $\mu_{R,\gamma}^{A}$, $\nu_R$ for $\nu_{R,\gamma}^{A}$, and $M_R$ for $M_{R,\gamma}^{A}$, unless otherwise stated.

\section{Equidistribution}
In this section, we establish the equidistribution of large horospheres.
The proof is adapted from that of \cite[Theorem~1.1]{Mir07a}
\begin{theorem} \label{thm:equidisOnP1M}
    We have weak convergence of probability measures on $\cP^1\cM_{g}$
    \[
        \frac{\nu_R}{M_R}
        \Rightarrow
        \frac{\nu_{g}}{b_{g}}
    \]
    as $R \to \infty$.
\end{theorem}
The following immediate corollary is exceedingly useful late on.
\begin{corollary} \label{cor:equidisOnM}
    We have weak convergence of probability measures on $\cM_{g}$
    \[
        \frac{\mu_R}{M_R}
        \Rightarrow
        \frac{B(X)}{b_{g}} \, \mu_{\WP}
    \]
    as $R \to \infty$.
\end{corollary}
\begin{proof}
    This follows from the fact that $\mu_R$ is the push-forward of $\nu_R$ via $\cP^1\cM_{g} \to \cM_{g}$ and Theorem~\ref{thm:equidisOnP1M}.
\end{proof}

The proof of Theorem~\ref{thm:equidisOnP1M} rests on the following series of propositions. 

Let $\nu$ be a weak limit of $(\nu_{R}/M_{R})_{R>0}$.
\begin{proposition} \label{prop:nuIsInvariant}
    The measure $\nu$ is invariant under the earthquake flow.
\end{proposition}

\begin{proposition} \label{prop:nuIsLebesgue}
    The measure $\nu$ is absolutely continuous with respect to $\nu_{g}$.
\end{proposition}

\begin{proposition} \label{prop:nuIsProba}
    The measure $\nu$ is a probability measure.
\end{proposition}

\begin{proof}[Proof of Theorem~\ref{thm:equidisOnP1M}]
    Proposition~\ref{prop:nuIsInvariant}, Proposition~\ref{prop:nuIsLebesgue}, and Theorem~\ref{thm:equidisOnP1M} imply that $\nu$ and $\nu_{g}$ differ by a multiplicative constant, and it follows from Proposition~\ref{prop:nuIsProba} that this constant is $1$.
\end{proof}

Proposition~\ref{prop:nuIsInvariant} is immediate.
\begin{proof}[Proof of Proposition~\ref{prop:nuIsInvariant}]
    This follows from the fact that $\nu_R$ is invariant under the earthquake flow (since $\nu_{g}$ is).
\end{proof}

For the rest of this section we shall prove Proposition~\ref{prop:nuIsLebesgue} and \ref{prop:nuIsProba}, which are more technical.

\subsection{Escape to infinity?}
In this subsection, we prove Proposition~\ref{prop:nuIsProba}.
The key ingredient is the following non-divergence result for the earthquake flow due to Y.~Minsky and B.~Weiss.
\begin{theorem}[{\cite[Theorem~E2]{MW02}, \cite[Corollary~5.12]{Mir07a}}] \label{thm:Mir07.5.11}
    For any $c>0$, there exists $\epsilon > 0$, depending only on $c$, such that for any $x \in\cT_{g}$ and any $\lambda\in\cML_{g}$, the following dichotomy holds:
    \begin{enumerate}
        \item[1.]
            There exists a simple closed curve $\alpha$ disjoint from $\lambda$, and $\ell_x(\alpha) < \epsilon$.

        \item[2.]
            We have
            \[
                \liminf_{T\to\infty} \frac{|\{ t \in [0,T] : \pi(\tw_\lambda^{\mkern1mu t}(x)) \in \cM_{g}^{\geq \epsilon} \}|}{T}
                >
                1-c
            \]
            where $\pi \colon \cT_{g} \to \cM_{g}$ is the natural projection, and $\cM_{g}^{\geq \epsilon}$ is the compact subset of $\cM_{g}$ consisting of all surfaces whose shortest closed geodesic has length at least $\epsilon$.
    \end{enumerate}
\end{theorem}

\begin{proof}[Proof of Proposition~\ref{prop:nuIsProba}]
    It is enough to prove that for any $\delta > 0$, we can find a compact subset $K_\delta$ of $\cP^1\cM_{g}$ such that
    \[
        \liminf_{R\to\infty} \frac{\nu_R(K_\delta)}{M_R}
        \geq
        1-\delta.
    \]
    The strategy is to show that there exists $\epsilon > 0$, depending only on $\delta$, such that the pre-image of $\cM_{g}^{\geq \epsilon}$ under $\cP^1\cM_{g} \to \cM_{g}$ possess the desired property.
    In other words,
    \[
        \liminf_{R\to\infty} \frac{\mu_R(\cM_{g}^{\geq \epsilon})}{M_R}
        \geq
        1-\delta.
    \]

    Taking $c = \delta/2$, Theorem~\ref{thm:Mir07.5.11} allows us to write $\tilde{\cS}_R \subset \cT_{g}$ as the disjoint union of $\tilde{\cS}_1$ and $\tilde{\cS}_2$ corresponding to the two possibilities. For convenience, we shall adapt the convention that $\bar{\cS}_{*}$ (resp.\ $\cS_{*}$) denotes the image of $\tilde{\cS}_{*}$ under $\cT_{g} \to \cM_{g}$ (resp.\ $\cT_{g} \to \cM_{g}^{\gamma}$), where $*$ is a certain index. 

    First, we show that $\mu_R(\bar{\cS}_1) \leq \mu_R(\cS_1) = \lo(M_R)$ as $R\to\infty$ even when $A = \varDelta^{k-1}$ (the subset of the simplex that we choose to define $\mu_R$ is the whole simplex).
    For any point in $\tilde{\cS}_1$, at least one of the following holds:
    \begin{enumerate}
        \item[1.1.]
            $\alpha$ is freely homotopic to $\gamma_i$ for some $1\leq i \leq k$.

        \item[1.2.]
            $\alpha$ is disjoint from $\gamma_1,\dots,\gamma_k$.
    \end{enumerate}
    Thus $\tilde{\cS}_1$ can be written as the union of $\tilde{\cS}_{1,1}$ and $\tilde{\cS}_{1,2}$ corresponding to the two cases above.
    To simplify the notation, in the following estimates of $\mu_R(\cS_{1,1})$ and $\mu_R(\cS_{1,2})$ we assume that $\gamma$ is primitive, i.e.\ $m_1 = \cdots = m_k = 1$ (the calculation differs from the general case only by a multiplicative constant).

    For each $i$, the corresponding horospherical volume of $\cS_{1,1}$ can be estimated by taking $f$ in Theorem~\ref{thm:MirIntFor} to be the indicator function
    \[
        \mathds{1} \{ (x_1,\dots,x_k) \in \RR_{\geq 0}^k : R \leq x_1 + \cdots + x_k \leq R + h,\ x_i < \epsilon \},
    \]
    and we obtain
    \[
        \mu_{R}({\cS}_{1,1})
        \leq
        \sum_{i=1}^{k} \lim_{h \to 0} \frac{\kappa_\gamma}{h} \int_0^\epsilon dx_i \int_{\varDelta_{[R-x_i, R+h-x_i]}^{k-1}} dx_1 \cdots dx_{i-1}dx_{i+1} \cdots dx_k \, P_\gamma(x_1,\dots,x_k),
    \]
    where $\varDelta_{[R-x_i, R+h-x_i]}^{k-1} \coloneqq \{ (x_1,\dots, x_{i-1}, x_{i+1}, \dots, x_k) \in \RR_{\geq 0}^{k-1} : R \leq x_1 + \cdots + x_k \leq R + h \}$.
    Since $P_\gamma$ is a polynomial of degree $2d - k$ and $x_1\cdots x_k$ is a factor of $P_\gamma$, we have $\mu_R(\cS_{1,1}) = \bO(\epsilon^2 R^{2d-3})$.

    We now suppose that $\alpha$ is disjoint from $\gamma_1,\dots,\gamma_k$.
    Denote by $(\gamma,\alpha)$ the ordered multicurve $(\gamma_1, \dots, \gamma_k, \alpha)$.
    Again, by applying Theorem~\ref{thm:MirIntFor}, $f$ being the indicator function
    \[
        \mathds{1} \{ (x_1,\dots,x_k, y) \in \RR_{>0}^{k+1} : R\leq x_1 + \cdots + x_k \leq R + h,\ y < \epsilon \},
    \]
    we obtain the corresponding horospherical volume
    \[
        \lim_{h \to 0} \frac{\kappa_\gamma}{h} \int_0^\epsilon dy \int_{\varDelta_{[R, R+h]}^{k}} dx_1 \cdots dx_k \, P_{(\gamma,\alpha)}(x_1,\dots,x_k,y)
    \]
    where $\varDelta_{[R,R+h]}^{k} \coloneqq \{ (x_1,\dots,x_k) \in \RR_{\geq 0}^k : R \leq x_1 + \cdots + x_k \leq R + h \}$.
    Since $P_{(\gamma,\alpha)}$ is a polynomial of degree $2d - k - 1$ of which $y$ is a factor, and there are only finitely many topological types of $(\gamma,\alpha)$, we have $\mu_R(\cS_{1,2}) = \bO(\epsilon^2 R^{2d-3})$.

    Using Corollary~\ref{cor:M}, we deduce
    \[
        \frac{\mu_R(\bar{\cS}_1)}{M_R}
        \leq
        \frac{\mu_R(\bar{\cS}_{1,1}) + \mu_R(\bar{\cS}_{1,2})}{M_R}
        \leq
        \frac{\mu_R(\cS_{1,1}) + \mu_R(\cS_{1,2})}{M_R}
        =
        \bO(\epsilon^2 R^{-2})
        =
        \lo(1)
    \]
    as $R \to \infty$.

    Let us now consider $\bar{\cS}_2$.
    One observes that every $p \in \bar{\cS}_2$ lies in a unique $1$-periodic earthquake flow orbit along the direction $\gamma_p \coloneqq \gamma_1/\ell_p(\gamma_1) + \cdots + \gamma_k/\ell_p(\gamma_k)$, and $\bar{\cS}_2$ can be written as the disjoint union of such orbits.
    (If one completes $\gamma$ to a pants decomposition by adding $d-k$ simple curves, such orbits are parallel straight lines in the $(\tau_{\gamma_1},\dots,\tau_{\gamma_k})$-coordinates plane.)
    By Theorem~\ref{thm:Mir07.5.11},
    \[
        |\{ t \in [0,1] : \pi(\tw_{\gamma_p}^{\mkern1mu t}(p)) \in \cM_{g}^{\geq \epsilon} \}|
        >
        1-\delta/2,
    \]
    for all $p \in \cS_2$.
    Thus,
    \[
        \frac{\mu_R(\bar{\cS}_2 \cap \cM_{g}^{\geq \epsilon})}{\mu_R(\bar{\cS}_2)}
        \geq
        1 - \delta/2.
    \]
    Therefore,
    \begin{align*}
        \frac{\mu_R(\cM_{g}^{\geq \epsilon})}{M_R}
        & =
        \frac{\mu_R(\cM_{g}^{\geq \epsilon} \cap \bar{\cS}_1) + \mu_R(\cM_{g}^{\geq \epsilon} \cap \bar{\cS}_2)}{M_R} \\
        & =
        0 + \frac{\mu_R(\cM_{g}^{\geq \epsilon} \cap \bar{\cS}_2)}{\mu_{R}(\bar{\cS}_2)} \frac{M_R-\mu_{R}(\bar{\cS}_1)}{M_R}
        \geq
        (1 - \delta/2)(1 - \lo(1))
    \end{align*}
    as $R \to \infty$. The $\lo(1)$ term can be made, e.g.\ smaller than $\delta/2$, by increasing $R$. The proof is thus complete. 
\end{proof}
\subsection{Absolute continuity}
In this subsection, we prove Proposition~\ref{prop:nuIsLebesgue}.

We use the following notation throughout the subsection.
Let $d$ denote $3g-3$.
We write $f = \bO_{K}(g)$ if there exists $C > 0$, depending only on $K$, such that $f \leq C g$, and we write $f = \Theta_K(g)$ if there exists $C$, depending only on $K$, such that $(1/C) g \leq f \leq C g$.

The key to the proof are the following estimates.
\begin{theorem}[\cite{Mir07a}, \cite{Mir}, \cite{AH21}] \label{thm:volB}
    Let $\epsilon \in (0,1)$, and let $K$ be a compact subset of $\cM_{g}$. We have
    \begin{enumerate}
        \item[1.]
            For any $x \in K$, $\mu_{\WP}(\BB_{x}(\epsilon)) = \Theta_K(\epsilon^{2d})$.

        \item[2.]
            For any $x \in \pi^{-1}(K) \subset \cT_{g}$, $\mu_R(\BB_x(\epsilon)) = \bO_K(\epsilon^{2d-1}/R)$.
    \end{enumerate}
\end{theorem}
\begin{remark} \label{rem:esti}
    The first part of the preceding theorem is \cite[Theorem~5.5.a]{Mir07a}.
    Mirzakhani proved the second part in the case when $\gamma$ is a simple closed curve \cite[Theorem~5.5.b]{Mir07a}, and claimed a more general version without proof \cite[Proposition~2.1.b]{Mir}.
    The proof of \cite[Theorem~5.5.b]{Mir07a} is concise and not easy to follow.
    See also the footnote on p.\ 390 in \cite{Wri20}.
    A much stronger estimate is obtained by Arana-Herrera in a different approach \cite[Proposition~1.5]{AH21}.
\end{remark}

The rest of the proof of Proposition~\ref{prop:nuIsLebesgue} can be adapted from Mirzakhani's original proof in the case when $\gamma$ is simple.
Let us sketch her arguments for the sake of self-containedness.
\begin{corollary} \label{cor:nuRnugn}
    Let $U\subset \PP(\cML_{g})$ be open, $K \subset \cT_{g}$ be compact, $x \in K$, and $p \colon \cT_{g} \times \PP(\cML_{g}) \to \cP^1\cM_{g}$ be the natural projection. For $\epsilon \in (0,1)$, we have
    \[
        \frac{\nu_R(p(\BB_x(\epsilon) \times U))}{M_R}
        =
        \bO_{K}( \nu_{g}(\BB_x(\epsilon) \times U_x) )
    \]
    where $U_x \coloneqq \{ \lambda \in \cML_{g} : \ell_x(\lambda) \leq 1,\ [\lambda] \in U \}$.
\end{corollary}
\begin{proof}
    It is enough to prove this for $A = \varDelta^{k-1}$.
    By \eqref{eq:1+-e}, for any $y \in \BB_x(\epsilon)$, we have
    \begin{equation} \label{eq:1+-eTh}
        (1 - \epsilon)^{2d} \cdot \mu_{\Th}(U_x)
        \leq
        \mu_{\Th}(U_y)
        \leq
        (1+\epsilon)^{2d} \cdot \mu_{\Th}(U_x),
    \end{equation}
    and so
    \begin{align*}
        & \#\{ \alpha \in \Mod_{g} \cdot \, \gamma : [\alpha] \in U,\ \ell_y(\alpha) = R \text{ for some }y\in\BB_x(\epsilon)\} \\
        & \leq \#\{ \alpha \in \Mod_{g} \cdot \, \gamma : [\alpha] \in U,\ (1-\epsilon) R \leq \ell_x(\alpha) \leq (1+\epsilon) R\} \\
        & = \bO_K(\epsilon R^{2d}\mu_{\Th}(U_x)).
    \end{align*}
    Hence Theorem~\ref{thm:volB}.2 implies that $\nu_R(p(\BB_x(\epsilon) \times U)) = \bO_K(\epsilon^{2d} R^{2d-1} \mu_{\Th}(U_x))$. The result now follows from Theorem~\ref{thm:volB}.1 and Corollary~\ref{cor:M}.
\end{proof}
We need one further technical lemma.
\begin{lemma} \label{lem:Mir07a.5.10}
    Let $K$ be a compact subset of $\cP^1\cT_{g}$. For any $N \subset K$ with $\nu_{g}(N) = 0$, and any $\epsilon > 0$, there exists an open cover $\{ \BB_{X_i}(r_i) \times U_i : i \in \ZZ_{\geq 1} \}$ of $N$, where for all $i$, $X_i \in \cT_{g}$, $r_i \in (0,1)$, and $U_i \subset \PP(\cML_{g})$ is open, such that
    \[
        \sum_{i\geq 1} \nu_{g}(\BB_{X_i}(r_i) \times U_i)
        \leq
        \epsilon.
    \]
\end{lemma}
\begin{proof}
    Fix a choice of Fenchel--Nielsen coordinates. There exists an open cover $\{ B_{X_i}(r_i) \times U_i : i \in \ZZ_{\geq 1} \}$ of $N$, where $B_{X_i}(r_i)$ is the Euclidean ball of radius $r_i$ centered at $X_i$, such that $\sum_{i \geq 1} \nu_g(B_{X_i}(r_i) \times U_i) \leq \epsilon$, and $\sup_{i\geq 1} r_i$ can be made as small as we please (since $\nu_g$ is a Lebesgue class measure).
    It follows from the compactness of $K \times [0,1]$ that there exists a constant $s$ depending only on $K$ such that $B_{x}(r) \subset \BB_x(s\cdot r)$ for any $x\in K$ and any $r \in [0,1]$.
    By Theorem~\ref{thm:volB}.1, there exists a constant $s'$ depending only on $K$ such that $\mu_\WP(\BB_{x}(s \cdot r)) \leq s' \cdot \mu_{\WP}(B_{x}(r))$ for any $x\in K$, and any $r < 1/2s$. Therefore, by \eqref{eq:1+-eTh}
    \[
        \sum_{i \geq 1} \nu_g(\BB_{X_i}(s \cdot r_i) \times U_i)
        \leq
        3^{2d} s'\sum_{i \geq 1}\nu_g(B_{X_i}(r_i) \times U_i)
        \leq
        3^{2d} s'\epsilon
        =
        \bO_K(\epsilon)
    \]
    and the lemma follows.
\end{proof}

\begin{proof}[Proof of Proposition~\ref{prop:nuIsLebesgue}]
    It is sufficient, as before, to consider the case in which $A$ is the whole simplex $\varDelta^{k-1}$.
    Let $N \subset \cP^1\cM_{g}$ with $\nu_{g}(N) = 0$. By Proposition~\ref{prop:nuIsProba}, we may assume that $N$ is contained in a compact set $K \subset \cP^1\cM_{g}$.
    Lemma~\ref{lem:Mir07a.5.10} implies for any $\epsilon > 0$, there exists an open cover $\{ \BB_{X_i}(r_i) \times U_i : i \in \ZZ_{\geq 1} \}$ of $N$, such that
    \[
        \sum_{i \geq 1} \nu_{g}(\BB_{X_i}(r_i) \times U_i)
        \leq
        \epsilon.
    \]
    Hence, it follows from Corollary~\ref{cor:nuRnugn} that 
    \[
        \sum_{i \geq 1} \frac{\nu_R(\BB_{X_i}(r_i) \times U_i)}{M_R}
        =
        \sum_{i \geq 1} \bO_K(\nu_{g}(\BB_{X_i}(r_i) \times U_{X_i}))
        \leq
        \bO_K(\epsilon).
    \]
    The proof is thus complete.
\end{proof}
\section{Counting}
The main result of this section is the following theorem which is a refined version of \cite[Theoreom~1.1]{Mir08b}.
\begin{theorem} \label{thm:countingCurves}
    Let $X\in \cM_{g}$, $\gamma = (m_1\gamma_1,\dots,m_k\gamma_k)$ be an ordered multicurve, and $A \subset \varDelta^{k-1}$ be open.
    We have
    \[
        \#\{ \alpha \in \Mod_{g} \cdot \,\gamma : \ell_X(\alpha) \leq R,\ \hat{\ell}_X(\alpha) \in A \}
        \sim
        C_{\gamma}^{A} \, \frac{[\Stab(\overline{\gamma}) : \Stab(\gamma)]}{2d} \frac{B(X)}{b_{g}} \, R^{2d}.
    \]
    as $R\to\infty$.
\end{theorem}
By virtue of Lemma~\ref{lem:bij}, this multicurves counting problem can be transformed to a counting problem on $\cM_{g}^{\gamma}$.
Let us begin by introducing some definitions that we need to state our counting result on $\cM_{g}^{\gamma}$.

The \emph{horoball} (on $\cM_{g}^{\gamma}$) is defined by
\[
    \cB_{R}
    \coloneqq
    \{ X \in \cM_{g}^{\gamma} : \ell_X(\gamma) \leq R,\ \hat{\ell}_X(\gamma) \in A \}
    =
    \bigcup_{0 < r \leq R} \cS_{r}
\]
and its associated measure $\mu_{\leq R}$ is defined by the formula
\[
    \mu_{\leq R}(U)
    \coloneqq
    \int_{0}^{R} \mu_{r}(U) \, dr
    =
    \mu_{\WP}(U \cap \cB_{R})
\]
where $U$ is any open subset of $\cM_{g}^{\gamma}$.
By abuse of notation, we shall also use $\mu_{\leq R}$ to denote the measure on $\cM_{g}$ defined by the formula
\[
    \nu_{\leq R}(U)
    \coloneqq
    \int_0^R \mu_r(U) \, dr
    =
    \mu_{\WP}(U \cap \pi_\gamma(\cB_R))
\]
for any open subset $U$ of $\cM_{g}$.
Let $X\in\cM_{g}$ and let $N(R)$ denote the number of pre-images of $X$ under $\pi_\gamma \colon \cM_{g}^\gamma \to \cM_{g}$ which lie within the horoball $\cB_R \subset \cM_{g}^\gamma$, i.e., $N(R) \coloneqq \#\{ \pi_\gamma^{-1}(X) \cap \cB_R \}$.

We have the following counting result on $\cM_{g}^{\gamma}$.
\begin{theorem} \label{thm:countingSurfaces}
    Let $X\in \cM_{g}$, $\gamma = (m_1\gamma_1,\dots,m_k\gamma_k)$ be an ordered multicurve, and $A \subset \varDelta^{k-1}$ be open. Then we have
    \[
        N(R)
        \sim
        C_{\gamma}^{A} \, \frac{[\Stab(\overline{\gamma}) : \Stab(\gamma)]}{2d} \frac{B(X)}{b_{g}} \, R^{2d}.
    \]
    as $R\to\infty$.
\end{theorem}
As an immediate corollary, we get the main result of this section:
\begin{proof}[Proof of Theorem~\ref{thm:countingCurves}]
    This follows at once from Theorem~\ref{thm:countingSurfaces} and Lemma~\ref{lem:bij}.
\end{proof}

We introduce a family of subsets $A_{a,b}$ of $\varDelta^{k-1}$, indexed by $a=(a_1,\dots,a_{k-1}) \in [0,1]^{k-1}$ and $b = (b_1,\dots,b_{k-1}) \in [0,1]^{k-1}$ such that $a_i < b_i$ for all $1\leq i \leq k-1$, and defined by
\[
    A_{a,b}
    \coloneqq
    \left\{ (x_1,\dots,x_k) \in \varDelta^{k-1} : a_i \leq x_i \leq b_i,\ \forall 1 \leq i \leq k-1 \right\}.
\]
To prove Theorem~\ref{thm:countingSurfaces}, it is enough to check the case when $A = A_{a,b}$ for all $a,b$.
In order to abbreviate our formulas, for the rest of this section we write
\[
    A
    \coloneqq
    A_{a,b},
    \qquad
    A_+
    \coloneqq
    A_{\frac{1-\epsilon}{1+\epsilon}a, \frac{1+\epsilon}{1-\epsilon}b},
    \qquad
    A_-
    \coloneqq
    A_{\frac{1+\epsilon}{1-\epsilon}a, \frac{1-\epsilon}{1+\epsilon}b},
\]
where we adopt the convention that $\frac{1+\epsilon}{1-\epsilon}b_i = 1$ if $\frac{1+\epsilon}{1-\epsilon} b_i > 1$,
and we write $\cB_{R}^{+} \coloneqq \cB_{R}^{A_+}$, $\mu_{\leq R}^{+} \coloneqq \mu_{\leq R}^{+}$, etc. 
The reason for the choice of $A_+$ and $A_-$ is the following elementary lemma.
\begin{lemma} \label{lem:-+}
    Choose $\epsilon \in (0,1)$ small enough to ensure that $A_-$ and $A_+$ are well-defined, and let $x,y \in \cM_{g}^\gamma$ with $d_\Th(x,y) \leq \epsilon$. We have
    \begin{enumerate}
        \item[(1)]
            If $x \in \cB_{(1-\epsilon)R}^{-}$, then $y \in \cB_{R}$.

        \item[(2)]
            If $x \in \cB_{R}$, then $y \in \cB_{(1+\epsilon) R}^{+}$.
    \end{enumerate}
\end{lemma}
\begin{proof}
    Suppose that $x \in \cB_{R}$. It follows from the inequality~\eqref{eq:1+-e} that
    \[
        \ell_y(\gamma)
        \leq
        (1+\epsilon) \ell_x(\gamma)
        \leq
        (1+\epsilon) R
    \]
    and
    \[
        \frac{1-\epsilon}{1+\epsilon}\, a_i
        \leq
        \frac{(1-\epsilon) \ell_x(m_i\gamma_i)}{(1+\epsilon) \ell_x(\gamma)}
        \leq
        \frac{\ell_y(m_i\gamma_i)}{\ell_y(\gamma)}
        \leq
        \frac{(1+\epsilon) \ell_x(m_i\gamma_i)}{(1-\epsilon) \ell_x(\gamma)}
        \leq
        \frac{1+\epsilon}{1-\epsilon}\, b_i
    \]
    which show that $y \in \cB_{(1+\epsilon)R}^{+}$.
    Part~(1) can be proved in a similar manner.
\end{proof}
Now we are ready to prove our main result of the section.
\begin{proof}[Proof of Theorem~\ref{thm:countingSurfaces}]
    We can choose $\epsilon \in (0,1)$ such that $\BB_{Y_1}(\epsilon) \cap \BB_{Y_2}(\epsilon) = \varnothing$ for any distinct pre-images $Y_1,Y_2$ of $X$ under $\pi_\gamma \colon \cM_{g}^\gamma \to \cM_{g}$.
    Let us write
    \[
        N_-(R)
        \coloneqq
        \#\{ Y \in \pi_\gamma^{-1}(X) \subset \cM_{g}^\gamma : \BB_Y(\epsilon) \subset \cB_R \},
    \]
    for set of all $Y \in \cM_{g}^\gamma$ such that $Y$ projects to $X$ and the Thurston distance ball of radius $\epsilon$ centered at $Y$ is entirely included within the horoball $\cB_R \subset \cM_{g}^\gamma$.
    Furthermore we write
    \[
        N_+(R)
        \coloneqq
        \#\{ Y \in \pi_\gamma^{-1}(X) \subset \cM_{g}^\gamma : \BB_Y(\epsilon) \cap \cB_R \neq \varnothing \}
    \]
    for the set of all $Y \in \cM_{g}^{\gamma}$ that project to $X$ such that $\BB_Y(\epsilon)$ intersects $\cB_R$.
    By definition,
    \[
        N_-(R)
        \leq
        N(R)
        \leq
        N_+(R).
    \]
    It follows from Lemma~\ref{lem:-+} that
    \begin{equation} \label{eq:N+}
        N_+(R) \cdot \mu_\WP(\BB_X(\epsilon))
        \leq
        \mu_\WP(\pi_\gamma^{-1}(\BB_X(\epsilon)) \cap \cB_{(1+\epsilon)R}^{+})
        =
        \mu_{\leq (1+\epsilon)R}^+(\pi_\gamma^{-1}(\BB_X(\epsilon)))
    \end{equation}
    and
    \begin{equation} \label{eq:N-}
        \mu_{\leq (1-\epsilon)R}^-(\pi_\gamma^{-1}(\BB_X(\epsilon)))
        =
        \mu_\WP(\pi_\gamma^{-1}(\BB_X(\epsilon)) \cap \cB_{(1-\epsilon)R}^{-})
        \leq
        N_-(R) \cdot \mu_\WP(\BB_X(\epsilon)).
    \end{equation}
    For any open subset $U\subset\cM_{g}$,
    \begin{equation} \label{eq:mu<=U}
        \mu_{\leq R}(\pi_\gamma^{-1}(U))
        =
        [\Stab(\overline{\gamma}) : \Stab(\gamma)] \cdot \mu_{\leq R}(U).
    \end{equation}
    We deduce from \eqref{eq:N+}, \eqref{eq:N-}, and \eqref{eq:mu<=U} that
    \[
        \mu_{\leq(1-\epsilon)R}^{-}(\BB_X(\epsilon))
        \leq
        \frac{N(R) \cdot \mu_\WP(\BB_X(\epsilon))}{[\Stab(\overline{\gamma}) : \Stab(\gamma)]}
        \leq
        \mu_{\leq(1+\epsilon)R}^{+}(\BB_X(\epsilon)).
    \]
    where $\BB_X(\epsilon) \subset \cM_{g}$.
    Hence,
    \begin{align}
        \lim_{R \to \infty} \frac{\mu_{(1+\epsilon)R}^{+}(\BB_X(\epsilon))}{R^{2d}}
        & =
        \lim_{R \to \infty} \frac{1}{R} \int_0^{(1+\epsilon)R} \frac{\mu_{t}^{+}(\BB_X(\epsilon))}{R^{2d-1}} \, dt \\
        & =
        C^{+} \cdot \lim_{R \to \infty} \frac{1}{R^{2d}} \int_0^{(1+\epsilon)L} t^{2d-1} \frac{\mu_t^+(\BB_X(\epsilon))}{C^{+} \cdot t^{2d-1}} \, dt \label{eq:limMu+B/L2d}
    \end{align}
    where $C^+ \coloneqq C_{\gamma}^{A_+}$ is given by \eqref{eq:M.=CR2d-1}.
    By Corollary~\ref{cor:equidisOnM} and \ref{cor:M},
    \[
        \frac{\mu_{t}^{+}(\BB_X(\epsilon))}{C^{+} \cdot t^{2d-1}}
        =
        \frac{1}{b_{g}} \int_{\BB_X(\epsilon)} B(Y)\, dY + \lo(1)
    \]
    as $t\to\infty$, where $\lo(1)$ is bounded and the constant depends only on $\gamma$ and $A^+$.
    Thus \eqref{eq:limMu+B/L2d} is equal to
    \[
        \frac{(1+\epsilon)^{2d}\, C^+}{2d\, b_{g}} \int_{\BB_X(\epsilon)} B(Y)\, dY.
    \]
    Therefore,
    \[
        \frac{\mu_{\WP}(\BB_X(\epsilon))}{[\Stab(\overline{\gamma}) : \Stab(\gamma)]} \, \limsup_{R \to \infty} \frac{N(R)}{R^{2d}}
        \leq
        \frac{(1+\epsilon)^{2d}\, C^+}{2d\, b_{g}} \int_{\BB_X(\epsilon)} B(Y) \, dY,
    \]
    and similarly,
    \[
        \frac{(1-\epsilon)^{2d}\, C^-}{2d\, b_{g}} \int_{\BB_X(\epsilon)} B(Y) \, dY
        \leq
        \frac{\mu_{\WP}(\BB_X(\epsilon))}{[\Stab(\overline{\gamma}) : \Stab(\gamma)]} \, \liminf_{R \to \infty} \frac{N(R)}{R^{2d}}
    \]
    where $C^- \coloneqq C_{\gamma}^{A_-}$.
    Taking $\epsilon \to 0$, we obtain
    \[
        \lim_{R \to \infty} \frac{N(R)}{R^{2d}}
        =
        C_{\gamma}^{A} \, \frac{[\Stab(\overline{\gamma}) : \Stab(\gamma)]}{2d} \frac{B(X)}{b_{g}}.
    \]
    This established the theorem.
\end{proof}

\section{Statistics}
\begin{proof}[Proof of Theorem~\ref{thm:fixedTopo}]
    Theorem~\ref{thm:countingCurves} implies
    \[
        \lim_{R \to \infty} \PP(\hat{\ell}_{X,R,\gamma} \in A)
        =
        \frac{C_{\gamma}^{A}}{C_{\gamma}^{\varDelta^{k-1}}}.
    \]
    The assertion now follows from Corollary~\ref{cor:M}.
\end{proof}
\begin{example}
    If $\gamma = (\gamma_1, \dots, \gamma_{3g-3})$ is a pants decomposition, then $g(v) = 0$, $n(v) = 3$, and $V_{g(v),n(v)} = 1$ for all $v\in V$. Thus $P_{\gamma}(x_1,\dots,x_{3g-3}) = x_1\cdots x_{3g-3}$, and Theorem~\ref{thm:fixedTopo} reduces to Theorem~\ref{thm:Mir.1.2}.
    \begin{figure}[h]
        \centering
        \begin{overpic}[width = 0.55\textwidth]{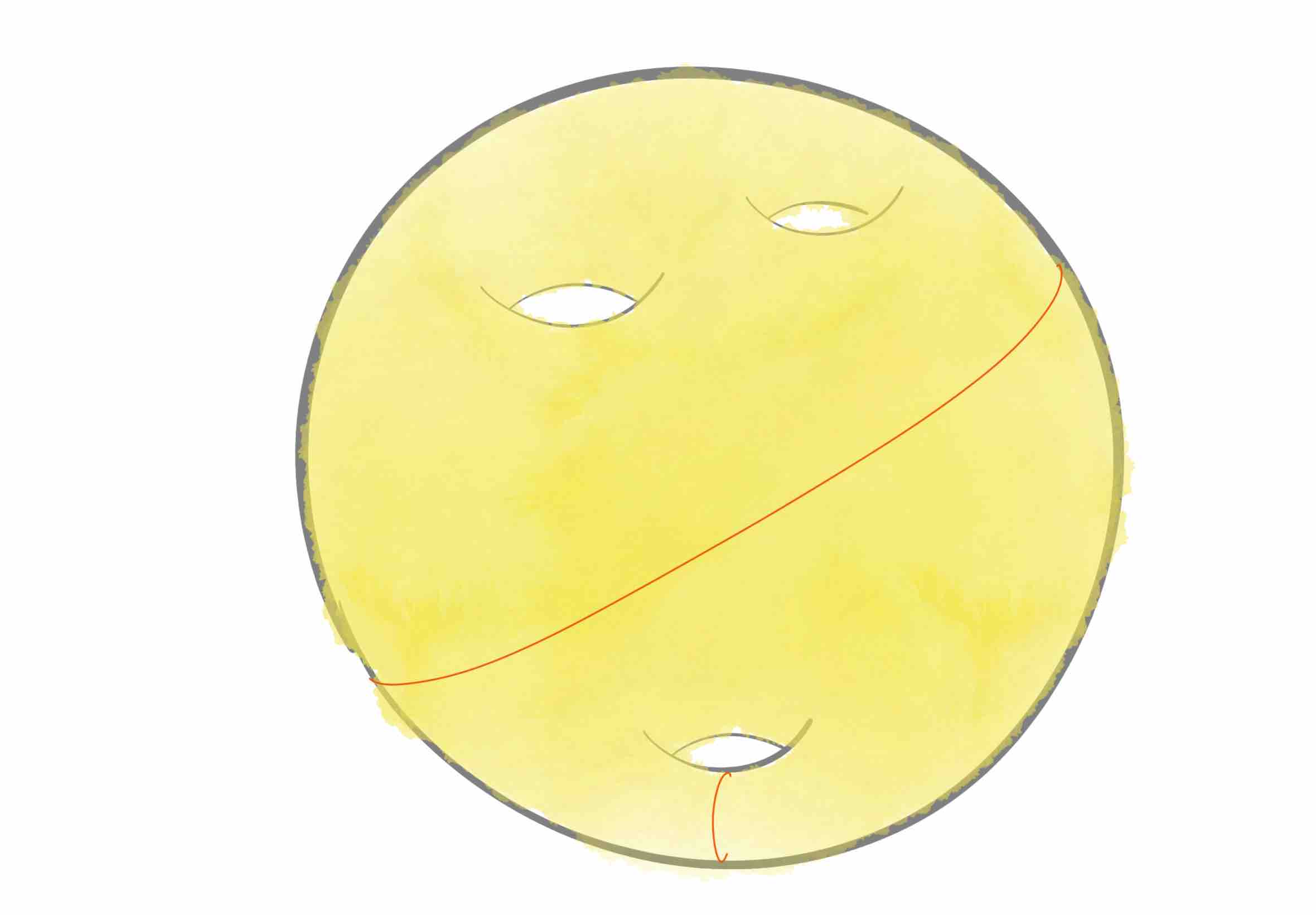}
            \put(54,0){$\gamma_1$}
            \put(70,33){$\gamma_2$}
        \end{overpic}
        \caption{Example~\ref{ex:2}} \label{fig:2}
    \end{figure}
\end{example}

\begin{example} \label{ex:2}
    Let $\gamma = (\gamma_1,\gamma_2)$, where $\gamma_2$ is separating and separates $\varSigma_{g}$ into a torus with a hole and a surface of type $(g-1,1)$, and $\gamma_1$ sits on the torus with a hole is non-separating as in Figure~\ref{fig:2}. Then its associated graph polynomial $\bar{P}_\gamma$ is equal to
    \[
        x_1 x_2 \cdot \bar{V}_{0,3}(x_1, x_1, x_2) \cdot \bar{V}_{g-1,1}(x_2)
        =
        \mathrm{constant} \cdot x_1 x_2^{6g-9}.
    \]
    This implies that in a random multi-geodesic of topological type $(\gamma_1,\gamma_2)$ on a hyperbolic surface of genus $g \gg 2$, the separating component is very likely to be much longer than the non-separating component.
\end{example}

\begin{example} \label{ex:3}
    Let $(\gamma_1, \gamma_2)$ be an ordered multicurve such that, for $i=1,2$, $\gamma_i$ is separating, and $\gamma_i$ bounds two surfaces of type $(g_i,1)$ (genus $g_i$ with $1$ boundary component) and $(g-g_1-g_2,2)$ respectively, as shown in Firgure~\ref{fig:3}. Then $\bar{P}_{\gamma}$ is
    \[
        x_1 x_2 \cdot \bar{V}_{g_1,1}(x_1) \, \bar{V}_{g_2,1}(x_2) \, \bar{V}_{g-g_1-g_2,2}(x_1,x_2)
        =
        \mathrm{constant} \cdot x_1^{6g_1-3} x_2^{6g_2-3} \cdot \bar{V}_{g-g_1-g_2,2}(x_1,x_2)
    \]
    where $\bar{V}_{g-g_1-g_2,2}$ is a symmetric polynomial.
    So in a typical multi-geodesic of type $(\gamma_1,\gamma_2)$, the first component is shorter than the second if $g_1 < g_2$.
    \begin{figure}[h]
        \centering
        \begin{overpic}[width = 0.55\textwidth]{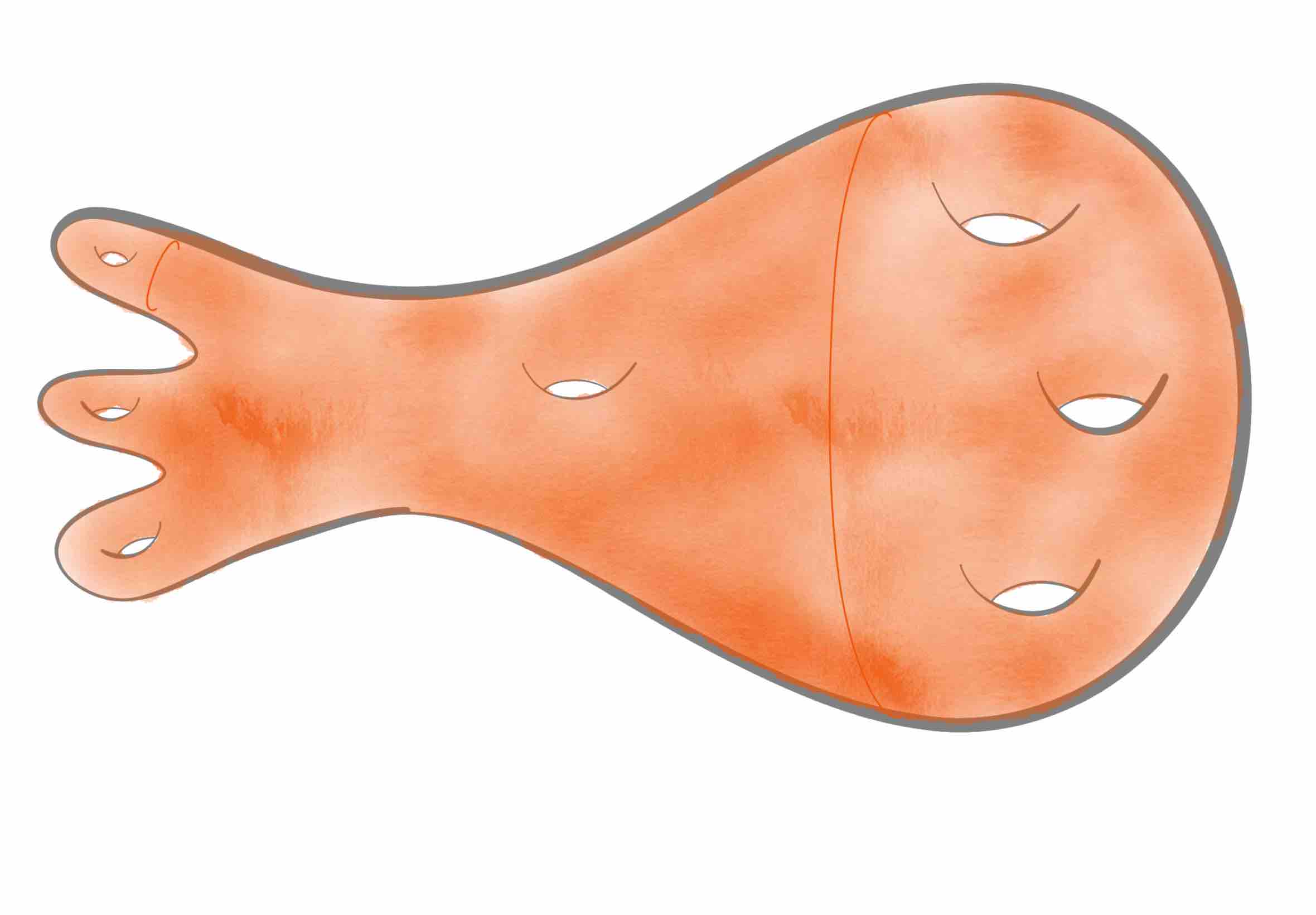}
            \put(12,54){$\gamma_1$}
            \put(58,28){$\gamma_2$}
            \put(1,55){$g_1$}
            \put(97,40){$g_2$}
        \end{overpic}
        \caption{Example~\ref{ex:3}} \label{fig:3}
    \end{figure}
\end{example}

\end{document}